\documentclass[preprints,article,accept,moreauthors,pdftex]{mdpi} 
\usepackage{smartdiagram}
\usepackage{graphicx}
\usepackage{verbatim}
\usepackage{marvosym}
\makeatletter
\tikzset{module/.style={%
		\pgfkeysvalueof{/smart diagram/module shape},
		thick,
		draw=\sm@core@bordercolor,
		top color=white,
		bottom color=\col,
		text=\sm@core@textcolor,
		minimum width=\sm@core@modulewidth,
		minimum height=\sm@core@moduleheight,
		font=\sm@core@modulefontsize,
		\sm@core@borderdecoration
	},
	diagram arrow type/.style={%
		\sm@core@arrowstyle,
		>=\sm@core@arrowtip,
		line width=\sm@core@arrowlinewidth,
		\col
	},%
}
\makeatother

\firstpage{1} 
\pubvolume{8}
\issuenum{8}
\articlenumber{888}
\pubyear{2021}
\copyrightyear{2021}
\datereceived{} 
\dateaccepted{} 
\datepublished{} 
\hreflink{https://doi.org/} 
\pdfoutput=1

\Title{Optimizing students' combinatorial-thinking skill through design-based research}

\Author{I. R. Ihsan $^{1}$\orcidA{}, and N. Karjanto $^{2,}$*\orcidB{}}

\AuthorNames{Iden Rainal Ihsan, and Natanael Karjanto}

\address{%
$^{1}$ \quad Departement of Mathematics Education, Faculty of Teacher Training and Educational Sciences, Universitas Samudra. Jalan Prof. Dr. Syarief Thayeb, Meurandeh, Langsa Lama, Langsa 24416, Aceh, Indonesia; irainalihsan@unsam.ac.id\\
$^{2}$ \quad Department of Mathematics, University College, Sungkyunkwan University, Natural Science Campus, 2066 Seobu-ro, Jangan-gu, Suwon 16419, Gyeonggi-do, Republic of Korea}

\corres{\Letter: \url{natanael@skku.edu} \hfill Updated \today}

\abstract{The aim of this study is to construct and compose an instructional design in combinatorial learning, particularly in the concept of counting. A composed design is expected to optimize students' combinatorial-thinking skill. This research was conducted at the first author’s workplace, in the topics of combinatorial and graph theory consisting of 12 students. This study consists of three phases: preliminary research, development, and assessment. In the first phase, we analyze the needs and context of learning activities, literature study, and develop a framework of thinking. In the second phase, we validates the instructional design for further improvement, classroom application, and improvement based on learning activities. In the third phase, we provide evaluation questions to students to find out the result of the application of the design related to the optimization of combinatorial-thinking. As a result, an instructional design in counting that begins with review of concept of finite sets and continues with small group discussion, class discussion, self-reflection, and evaluation.}

\begin{document}

\section{Introduction}

As prospective mathematics educators, students majoring in Mathematics Education must not only be able to master most of the concepts that will be taught later when they become teachers themselves but also need to possess a capability in teaching higher-order thinking skills~\cite{ref-ihsankosasih2021penelitian}. By simultaneously developing pedagogical skills, teachers will construct the well-known pedagogical content knowledge (PCK) as practicioners. Furthermore, together with beliefs, this PCK constitutes essential factors that will influence teaching practice in the classroom.
	
Combinatorics, one of the compulsory courses for them, contains rich topics that can be utilized to inculcate these skills. Combinatorics is an area of mathematics primarily deals with enumeration or counting of specified structures and particular properties of such structures. It also covers the existence of and the construction of those structures. The course covers abundant problem materials and it provides beneficial foundation for other interdisciplinary studies, including Natural and Computer Sciences, Engineering, a well as Probability and Statistics.~\cite{ref-lockwood2013model}

Future mathematics teachers ought to possess combinatorial reasoning and creativity skill in order to cultivate mathematics among their pupils. According to Rezaie and Gooya, there are four levels of understanding combinatorial-thinking~\cite{ref-rezaie2011what}. The first level is ``investigating some cases''. In this level, students can investigate certain and special cases of combinatorics. The second category is ``How am I sure that I have counted all the cases?''. In this level, students argue about their investigation. The third level is ``systematically generating all cases''. In this level, students can make a generalization problem systematically, from particular simple cases to more general complex. Finally, the fourth level is ``changing the problem into another combinatorial problem''. In this level, students can solve other problems related to and relevant to the previously worked problems. 

This study aims to construct and compose an instructional design in combinatorial learning, especially in the concept of counting. We expect to optimize students' combinatorial-thinking skill by the composed learning design. Furthermore, we also hope that students can construct comprehensive understanding during this learning process. It is expected that when they will become teachers themselves, they will not only set an an excellent example in teaching but also construct their students' understanding as well. Indeed, as emphasized by Pradipta et al., exemplary educators can select valid and practical teaching methods that emphasize comprehensive student understanding, and do not merely perform a set of exemplifications~\cite{ref-pradipta2021exem}. A differentiation in instructions could help these mathematics educators to make informed decisions about how to support students who bring various resources into the classroom~\cite{ref-tillema2021student}.

In an effort to compose a learning design that suits our goal, we adopt a methodological approach called the design-based research (DBR). DBR is an increasingly popular practical research method that assists both researchers and practicioners {to come together} in order to generate effective educational interventions and methods~\cite{ref-anderson2012design,ref-minichiello2021narrative}. DBR is an approach to educational research wherein creative and generative processes of design are intimately intertwined with theresearch process: the design is research-based and the research is design-based~\cite{ref-bakker2019design}. Related to that intertwining, Haagen-Schützenhöfer and Hopf argued that there are noteworthy contributions to the field of education~\cite{ref-haagens2020design}.

We choose DBR because it could serve as an effective bridge between research and practice in the sector of formal education~\cite{ref-anderson2012design,ref-vandermerwe2019design}. Using DBR, we are able to develop ready-to-use curriculum materials~\cite{ref-vandermerwe2019design} and it focuses on real-world problems~\cite{ref-zhao2021innovative,ref-dolmans2019how}. Since our objective is to optimize students' combinatorial-thinking, DBR is attractive since it focuses on optimizing teachers' interventions~\cite{ref-ivens2020does}. DBR can be designed by and for educators who seek to increase the impact, transfer, and translation of education research into improved practice. Additionally, since the theory of DBR could aid us in explaining of how complex phenomena in problem-based learning (PBL) interact and thereby add to our understanding of these phenomena, we also include PBL materials to our design~\cite{ref-santos2010using}.
\begin{figure}[H]
\begin{center}
\includegraphics[width=12 cm]{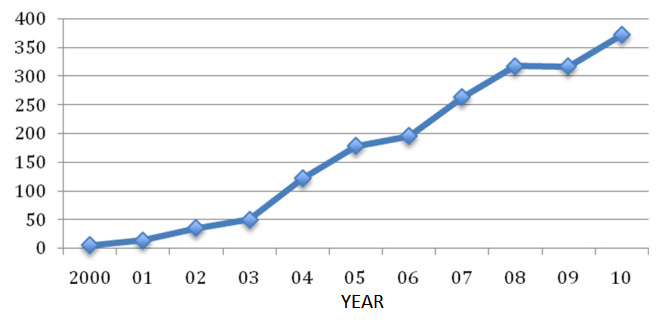}
\end{center}
\caption{A graph showing the number of articles covering or discussing DBR.} \label{fig1}
\end{figure}

Historically, DBR evolved near the beginning of the 21st century~\cite{ref-anderson2012design}. Two factors have promoted its emergence. First, DBR concerns about the low impact of educational research on educational practice. Second, the extended use of constructivist situated theories~\cite{ref-brown1992design}. For the past two decades, the growth of research on DBR has been steadily increasing, as can be seen in Figure~\ref{fig1}. There are several other different terms referring to DBR. In the Netherlands, it is known as ``developmental or development research'', which is put forward by researchers who focus on curriculum development, e.g., Freudenthal, Gravemeijer, Linjse, Romberg, and Van der Akker~\cite{ref-minichiello2021narrative}. Unlike in the United States, for which the term ``design experiments'' is more commonly used, the term DBR comes from researchers in the field of Cognitive Psychology~\cite{ref-brown1992design,ref-collins1992toward,ref-dbrc2003design}, and only recently that it has grown into its new name~\cite{ref-bakker2019design}. A DBR collective team has committed to theorizing and practicing design-based research in education emerged  under the name ``design-based research'' in 2003~\cite{ref-dbrc2003design}. It has since become a widely used umbrella term to designate research approaches that blend the design of educational interventions with educational and learning research~\cite{ref-minichiello2021narrative}. 

DBR is frequently described in the literature as being pragmatic, interventionist, and collaborative~\cite{ref-ford2017using}. Researchers of DBR reveal better educational solutions through a continuous improvement~\cite{ref-chen2020redeveloping}. Although DBR can enrich experimental and quasi-experimental research, the objectives differ in some important aspects. DBR is more than just proving that one pedagogical approach is more effective than another. The solution resulted from implementing DBR is in accordance with theoretical principles and knowledge of recent technological availability~\cite{ref-ivens2020does}. 

There is another view related to DBR. Other than those previously described, DBR can be also characterized as an approach which blends empirical educational research with the theory-driven design of learning environments~\cite{ref-brown1992design}. DBR output contributes with both theoretical knowledge and societal education~\cite{ref-zhao2021innovative}. DBR approach has a cyclic character in which thought experiments and teaching experiments alternate. A cycle consists of three phases, i.e., the preliminary, the teaching experiment phase or developing, and the phase of analysis or assessment \cite{ref-doorman2018design,ref-plopm2013educational}.

Based on the aim of this study, this paper is organized as follows. Section~\ref{methodology} covers the description of the methodology. We also provide details of DBR instructional design. Section~\ref{result} gives the results of our design implementation, where the level of students' combinatorial-thinking is displayed. Sections~\ref{discuss} and~\ref{conclude} provides discussion and conclusion, respectively.

\section{Materials and Methods}	\label{methodology}

\begin{figure}[h]
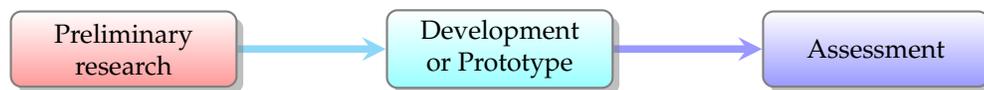

\begin{center}
\smartdiagramset{module minimum height=1cm,
	module minimum width=3cm,
	module x sep=5,back arrow disabled=true}
\smartdiagram[flow diagram:horizontal]{Preliminary\\ research,Development \\or Prototype,Assessment}
\end{center}
\caption{Three phases of DBR showing preliminary research, development, and assessment.\label{fig2}}
\end{figure}

We adopted DBR methodology. It consists of three stages: preliminary research, development, and assessment, as shown in Figure~\ref{fig2}~\cite{ref-plopm2013educational}.

During the preliminary research stage, several steps were involved. We first analyze the curriculum prerequisite for combinatorial-thinking skill and students' learning needs. The latter was obtained from the pre-test of combinatorial-thinking skill as well as from the previous learning outcomes on addition and multiplication rules. We then conduct a review on literature study on permutation and combination, e.g.,~\cite{ref-rezaie2011what,ref-herman2003counting,ref-bryant2003aspect}. The final step in this stage is designing learning framework which is arranged according to the task-based learning and guided-learning concepts~\cite{ref-willis1996framework,ref-ihsan2015penelitian}.

The second stage of development is also known as prototype. In this stage, we provide learning design to expert team for further validation. Based on validation results, we implement a revised learning design in class. Learning activity results were then made as revised learning design materials. At the assessment stage, we provide test questions to be completed by each student. The purpose is to determine the students' level of mastery in combinatorial-thinking.

We conducted this research at the first author's institution during the Spring semester of the 2017/2018 academic year. The selected course was called \emph{Combinatorics and Graph} (course code MT406619). The student body consists of twelve students, where five of them were males and the other seven were females. All of them have fulfilled the prerequisites courses for MT406619, i.e., Basic Mathematics, Logic and Sets, Number Theory, and Mathematical Thinking.

\section{Results}	\label{result}

We obtain two general outcomes. The first is a design of learning material in \emph{Permutation} that relates to combinatorial-thinking ability. The second is the evaluation result of the combinatorial-thinking ability level. The learning design in \emph{Permutation} consists of three main activities, i.e., initial, core, and closing. The initial activity of the lecture begins with the preparation of learning activities that consists of material review on general prerequisite concepts and in particular, directing students to review \emph{Finite Sets} themselves.
\begin{specialtable}[H] 
\caption{The result of combinatorial-thinking mastery level for each student.} \label{tab1}
\begin{center}
\begin{tabular}{ccccccccccccc}
\toprule
\textbf{Student}	& \textbf{1}	& \textbf{2}	& \textbf{3}	& \textbf{4}	& \textbf{5}	& \textbf{6}	& \textbf{7}	& \textbf{8}	& \textbf{9}    &\textbf{10}	& \textbf{11} & \textbf{12}\\
\midrule
Level of mastery in \\ combinatorial-thinking		& 4			& 4 & 3			& 4 & 4			& 4 & 2			& 3 & 4			& 4 & 4			& 4\\
\bottomrule
\end{tabular}
\end{center}
\end{specialtable}

During the core activities, students are directed to complete tasks in groups, where one group consists of two students. The given task is in the form of problems related to \emph{Permutation}. Problems are presented in three consecutive questions. The first question is related to a simple problem, assuming that each student can complete briefly. This inquiry accommodates the assessment of mastery Level~1 in students' combinatorial-thinking skill. The second question directs students to provide arguments for their answers to the first question. This second inquiry accommodates the assessment of mastery Level~2 in combinatorial-thinking skill. The third question directs students to summarize their findings, and thus combinatorial-thinking skill Level~3 can be assessed and observed.
\begin{figure}[H]
\begin{center}
\includegraphics[width=12 cm]{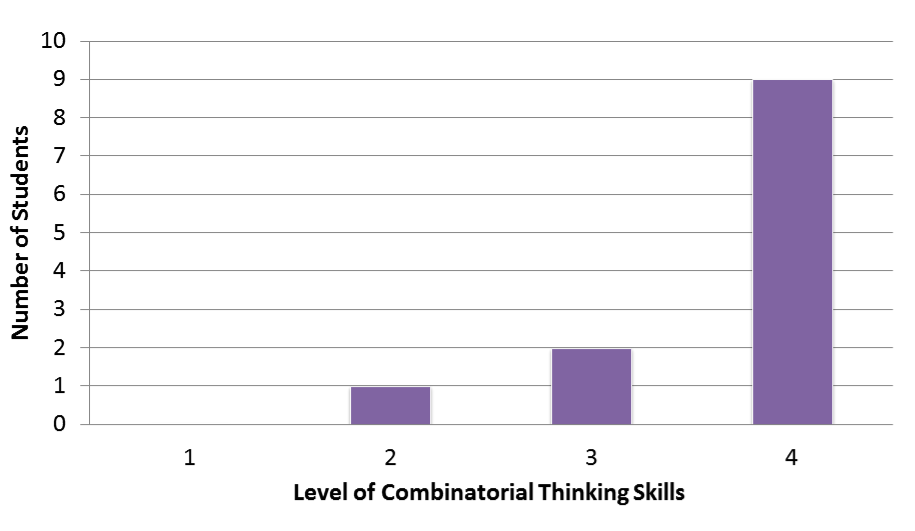}
\end{center}
\caption{A histogram summarizing the result of students' combinatorial-thinking skill.} 	\label{fig3}
\end{figure}  

In an effort to direct students in finding valid \emph{Permutation} concepts, we continue the learning process to the stage of class discussion. Two selected groups would present the results of their respective group discussions. The groups that do not present are invited to provide feedback. In the closing activity, each student is given a questionnaire. Each student has 15 minutes to complete it. The purpose of this questionnaire is to investigate students' mastery in \emph{Permutation}. To measure the mastery of students' combinatorial-thinking skill, particularly Level 4, we provide additional tasks that need to be completed outside regular contact hours. The assignments also include reviewing learning materials for the following in-class meeting. By analyzing the learning and evaluation activities, we obtain a profile of mastery in students' combinatorial-thinking skill, as shown in Table~\ref{tab1} and Figure~\ref{fig3}. 

\section{Discussion} \label{discuss}

Our study is in accordance to the recent work that elaborates on supporting new teachers~\cite{ref-brown2020supporting}. Indeed, nowadays teaching is a complex context and requires a variety of competencies. With this combinatorial-thinking ability, we can glean two lessons that might be useful for prospective-teacher students. The first is learning to solve complex problems, which one can start from simple examples and gradually progresses to more complex issues. Second, students can learn instructional designs that can be utilized as an alternative pedagogy in teaching when they become teachers themselves.

From the results in implementing learning design, we observed that the majority of our students have succeeded in reaching Level 4 of combinatorial-thinking skill. A sequence of discussions, for which students have opportunities to share experiences and expand their knowledge, contribute to this level of success. These discussion activities are certainly not limited to listening and following presentations. In addition to attempting problems and observing from one another, students do not only gain meaningful experience and grow in mathematical maturity (the ``making process'') but also to reflect on their own combinatorial-thinking skill (the ``sharing and openness process'')~\cite{ref-kim2020development}.

\section{Conclusion} \label{conclude}

We have implemented DBR in \emph{Combinatorics and Graph} course, and in particular designing learning material for the \emph{Permutation} topic. Using DBR, the learning process follows three stages of activity: initial, core, and closing. After completing any preparation, the initial activity proceeds with reviewing prerequisite concepts, which include \emph{Finite Sets}. The core activity consists of group and class discussions. The final learning activity comprises evaluation and out-of-class assignments.

Regarding the level of mastery in students' combinatorial-thinking skill, it turns out that not everyone reaches the optimum level of mastery. Nonetheless, the majority of them did reach Level 4, and out of 12 students, only one and two students attain Levels~1 and~2, respectively. 

\vspace{6pt} 

\authorcontributions{Conceptualization, I.R.I. and N.K.; methodology, I.R.I. and N.K. ; validation, N.K. ; formal analysis, N.K.; investigation, I.R.I ; resources, I.R.I; writing---original draft preparation, I.R.I; writing---review and editing, N.K.; supervision, N.K.; project administration, I.R.I;  All authors have read and agreed to the published version of the manuscript.}

\funding{This research was funded by the Indonesian Ministry of Research, Technology, and Higher Education, under the grant number 0045/E3/LL/2018.}

\informedconsent{Informed consent was obtained from all subjects involved in the study.}
 
\acknowledgments{On this occasion, the first author would like to thank the Directorate General of Research and Development Strengthening, Ministry of Research, Technology and Higher Education, Republic of Indonesia, for providing funding assistance to participate in the 9th International Congress on Industrial and Applied Mathematics 2019 through BSLN program.}

\conflictsofinterest{The authors declare no conflict of interest. The funders had no role in the design of the study; in the collection, analyses, or interpretation of data; in the writing of the manuscript, or in the decision to publish the results.} 

\end{paracol}

\reftitle{References}



\end{document}